\documentclass{amsart}
\usepackage{ifthen}
\usepackage{enumerate}

\numberwithin{equation}{section}
\def\epsilon{\varepsilon}
\def\parti#1#2{\frac{\partial #1 } {\partial #2} }

\def\de{\delta}
\def\ti{\tilde}

\def\ga{\gamma}
\def\uph{{ {}^h } }
\def\up#1{{ {}^{#1} } }

\def\beq{\begin{eqnarray}}
  \def\eeq{\end{eqnarray}}
\def\be{\beta}
\def\al{\alpha}
\def\ep{\epsilon}
\def\parttt{\frac{\partial }{\partial \tilde t} }
\def\partt{\frac{\partial }{\partial t} }
\def\parts{\frac{\partial }{\partial s } }
\def\phi{\varphi}
\def\R{\mathbb R}
\def\boundary{\partial}
\def\N{\mathbb N}
\def\M#1#2{\mathcal M^{#1}\left(#2\right)}

\def\H{\mathbb H}
\def\delt{\frac{\partial}{\partial t}}

\def\gradh{{{}^{{}^h}\!\nabla}}

\def\ep{\epsilon}

\def\part{\partial}

\def\curlL{\mathcal L}
\def\grad{\nabla}

\DeclareMathOperator{\vol}{vol}

\DeclareMathOperator{\Ric}{Ric}
\DeclareMathOperator{\Ricci}{Ric}

\DeclareMathOperator{\id}{id}
\def\loc{\text{\it{loc}}}
\def\counterword#1{%
  \ifthenelse{\ref{#1}=1}{one}{}%
  \ifthenelse{\ref{#1}=2}{two}{}%
  \ifthenelse{\ref{#1}=3}{three}{}%
  \ifthenelse{\ref{#1}=4}{four}{}%
  \ifthenelse{\ref{#1}=5}{five}{}%
  \ifthenelse{\ref{#1}=6}{sic}{}%
}

\newtheorem{theorem}{Theorem}[section]
\newtheorem{lemma}[theorem]{Lemma}
\newtheorem{corollary}[theorem]{Corollary}

\newtheorem{proposition}[theorem]{Proposition}
\theoremstyle{definition}
\newtheorem{definition}[theorem]{Definition}

\theoremstyle{remark}
\newtheorem{remark}[theorem]{Remark}
\def\fuhome{@math.fu-berlin.de}

\begin{document}
\title{Stability of Hyperbolic space under Ricci Flow}


\author{Oliver C. Schn\"urer}
\address{Oliver Schn\"urer:
  Fachbereich Mathematik und Statistik, 
  Universit\"at Konstanz,  
  78457 Konstanz, Germany}
\curraddr{}
\def\AmSeeHome{@uni-konstanz.de}
\email{Oliver.Schnuerer\AmSeeHome}


\author{Felix Schulze}
\address{Felix Schulze: 
  Freie Universit\"at Berlin, Arnimallee 6, 
  14195 Berlin, Germany}
\curraddr{}
\email{Felix.Schulze\fuhome}

\author{Miles Simon}
\address{Miles Simon: Universit\"at Freiburg, Eckerstra{\ss}e 1,
  79104 Freiburg i. Br., Germany}
\curraddr{}
\email{msimon@gmx.de}

\date{March 2010}

\dedicatory{}


\begin{abstract}
  We study the Ricci flow of initial metrics which are
  $C^0$-per\-tur\-ba\-tions of the {\it hyperbolic metric} on $\H^n$.  If
  the perturbation is bounded in the $L^2$-sense, and small enough in
  the $C^0$-sense, then we show the following: In dimensions four and
  higher, the scaled Ricci harmonic map heat flow of such a metric
  converges smoothly, uniformly and exponentially fast in all
  $C^k$-norms and in the $L^2$-norm to the hyperbolic metric as time
  approaches infinity. We also prove a related result for the Ricci
  flow and for the two-dimensional conformal Ricci flow.
\end{abstract}

\maketitle
\section{Introduction}
We investigate stability of hyperbolic space under Ricci flow
\begin{equation}
  \label{eq:RF}
  \begin{cases}
    \delt g_{ij}=-2\Ric(g(t)) &\text{on }\H^n\times(0,\infty),\\
    g(0)=g_0&\text{on }\H^n.
  \end{cases} 
\end{equation}
As hyperbolic space expands under Ricci flow, it is convenient to
consider the following modified Ricci flow
\begin{equation}
  \label{eq:RFresc}
  \begin{cases}
    \delt g_{ij}=-2\Ric(g(t))-2(n-1)g_{ij}(t) &\text{on
    }\H^n\times(0,\infty),\\ 
    g(0)=g_0&\text{on }\H^n.
  \end{cases}
\end{equation}
The hyperbolic metric $h$ of sectional curvature $-1$ is a stationary
point to \eqref{eq:RFresc}. 

Note that by Lemma \ref{lem:equiv}, up to rescaling, this flow
equation is equivalent to \eqref{eq:RF}. As \eqref{eq:RF} and
\eqref{eq:RFresc} are degenerate parabolic equations, we consider the
following modified (or rescaled) Ricci harmonic map heat flow which is
similar to DeTurck flow (\cite{DeTurck})
\begin{equation}
  \label{eq:RTFR}
  \begin{cases}
    \delt g_{ij}=-2\Ric(g(t))+\nabla_iV_j+\nabla_jV_i-2(n-1)g_{ij}(t)
    &\text{on }\H^n\times(0,\infty),\\
    g(0)=g_0&\text{on }\H^n,
  \end{cases}
\end{equation}
where $V_i=g_{ik}\left(\up g\Gamma^k_{rs}-\up h\Gamma
  ^k_{rs}\right)g^{rs}$ and $h$ is the hyperbolic metric on $\H^n$ of
sectional curvature equal to $-1$. Once again, up to rescaling, this
is equivalent to standard Ricci harmonic map heat flow. 

We consider perturbations that are close to hyperbolic space.
\begin{definition}
  \label{close def}
  Let $g$ be a metric on $\H^n$. Let $\epsilon>0$. Then $g$ is
  $\epsilon$-close to $h$ if 
  \[(1+\epsilon)^{-1}h\le g\le(1+\epsilon)h\] 
  in the sense of matrices. 
\end{definition}
Let $(\lambda_i)_{1\le i\le n}$ denote the eigenvalues of $(g_{ij})$
with respect to $(h_{ij})$. Then this is equivalent to
$(1+\epsilon)^{-1}\le\lambda_i\le1+\epsilon$ for $1\le i\le n$.

We denote with $\mathcal M^k(\H^n,I)$ the space of families
$(g(t))_{t\in I}$ of sections in the space of Riemannian metrics on
$\H^n$ which are $C^k$ on $\H^n\times I$. Similarly, we define
$\mathcal M^\infty$, $\mathcal M_{\loc}^k$ and use $\mathcal
M^k(\H^n)$ if the metric does not depend on $t$. We wish to point out
that we use $C^k$ on non-compact sets to denote the space, where
covariant derivatives with respect to the hyperbolic metric $h$ of
order up to $k$ are continuous and in $L^\infty$. We also use
$C^k_{\text{\it{loc}}}$. For our convenience, we define $\mathcal
M^\infty_0(\H^n,[0,\infty))$ to be the set of all metrics in $\mathcal
M^0(\H^n,(0,\infty)) \cap
M^0_{\loc}(\H^n,[0,\infty)) $ which are smooth for positive times and, when
restricted to time intervals of the form $[\delta,\infty)$,
$\delta>0$, are uniformly bounded in $C^k$ for any $k\in\N$.

We use $|Z|\equiv\up h|Z|$ to denote the norm of a tensor $Z$ with
respect to the hyperbolic metric $h$. Unless stated otherwise,
$B_R(0)$ denotes a geodesic ball around a fixed point in hyperbolic
space which we denote by $0$. $\Vert\cdot\Vert_{L^2}$ denotes the
$L^2$-Norm with respect to the hyperbolic metric $h$. Sometimes, we
write $x\to\infty$ instead of $|x|\to\infty$. Moreover, we use the
Einstein summation convention and denote generic constants by $c$.

Our main theorem is 
\begin{theorem}\label{thm:main}
  Let $n\ge4$. For all $K>0$ there exists
  $\epsilon_1=\epsilon_1(n,K)>0$ such that the following holds. Let
  $g_0\in\mathcal M^0(\H^n)$ satisfy $\int\limits_{\H^n}|g_0-h|^2 \,
  d\!\vol_h \leq K$ and $\sup\limits_{\H^n}|g_0-h| \leq
  \epsilon_1$. Then there exists a solution $g \in
  \mathcal{M}^\infty_0(\H^n, [0,\infty))$ to \eqref{eq:RTFR} such that
  \begin{equation*}
    \sup_{\H^n} |g(t)-h| \leq C(n, K)\cdot e^{-\frac 1{4(n+2)}t}.
  \end{equation*}
  Moreover, $g(t)\to h$ exponentially in $C^k$ as $t\to\infty$ for all
  $k\in\N$. 
\end{theorem}

There might be different solutions to the ones obtained by our
construction. The following theorem does not assume that the solution
in question comes from the theorem above.
\begin{theorem}\label{thm:sol given}
  Let $n\ge4$. For all $K>0$ there exists
  $\epsilon_1=\epsilon_1(n,K)>0$ such that the following holds. Let
  $g\in\mathcal M^\infty_0(\H^n\times[0,\infty))$ be a solution to
  \eqref{eq:RTFR} that satisfies $\int\limits_{\H^n}|g(0)-h|^2 \,
  d\!\vol_h \leq K$, $\lim\limits_{r\to\infty}\Vert
  g(0)-h\Vert_{L^\infty(\H^n\setminus B_r(0))}=0$ and
  $\sup\limits_{\H^n}|g(t)-h| \leq \epsilon_1$ for all $t\ge0$. Then
  \begin{equation*}
    \sup_{\H^n} |g(t)-h| \leq C(n, K)\cdot e^{-\frac 1{4(n+2)}t}.
  \end{equation*}
  Moreover, $g(t)\to h$ exponentially in $C^k$ as $t\to\infty$ for all
  $k\in\N$. 
\end{theorem}

If $g(0)\in\mathcal M^\infty_{\loc}\left(\H^n\right)$, solutions of
\eqref{eq:RTFR} correspond to solutions of \eqref{eq:RFresc}
\begin{theorem}\label{thm:Ricci}
  Let $n\ge4$. Let $g$ be a solution to \eqref{eq:RTFR} as in Theorem
  \ref{thm:main} or Theorem \ref{thm:sol given}. Assume in addition
  that $g$ is smooth. Then there exists a smooth family of
  diffeomorphisms of $\H^n$, $\phi_0=\rm{id}_{\H^n}$, such that for
  $\tilde{g}(t):=\phi_t^*g(t)$ the family $(\tilde{g}(t))_{t\geq 0}$
  is a smooth solution to \eqref{eq:RFresc} satisfying
  $$\tilde{g}(t)\to(\phi_\infty)^*h\quad\text{in } 
  \mathcal M^\infty\big(\H^n\big)\text{ as }t\to\infty$$ for some
  smooth diffeomorphism $\phi_\infty$ of $\H^n$ which satisfies
  $\phi_t\to\phi_\infty$ in $C^\infty\left(\H^n,\H^n\right)$ as
  $t\to\infty$ and, if $\lim\limits_{r\to\infty}\Vert
  g(0)-h\Vert_{L^\infty\left(\H^n\setminus B_r(0)\right)}=0$,
  $$|\phi_\infty(x)-x|\to0\quad\text{as }|x|\to\infty.$$
\end{theorem}

\begin{remark}
  All of the above results carry over directly if one replaces
  $(\mathbb{H}^n,h)$ by any complete Riemannian manifold $(M,h)$, where
  $h$ has sectional curvature equal to $-1$. 
\end{remark}

Linearised stability of hyperbolic space under Ricci flow has been
investigated before by V. Suneeta \cite{Suneeta}.  R. Ye considered
stability of negatively curved manifolds on compact spaces in the
paper \cite{YeEinstein}. H. Li and H. Yin \cite{LiYinHypRic} have
shown stability of hyperbolic space in dimensions $n \geq 6$ under the
assumptions that the deviation of the curvature of the initial metric
from hyperbolic space decays exponentially {\it and} the initial
metric is close to $h$ (in the sense of Definition \ref{close def}).

Similar results and methods to those found in this paper may be found
in the authors' paper \cite{Riccistab} addressing the stability of
Euclidean space under Ricci flow. For further references, we refer to
the introduction therein.

Here we  outline the proofs of the main results of this paper, and
explain the structure of these proofs and of the paper.

In the first part of the paper (chapters \counterword{exist sec} and
\counterword{conv sec}) we consider the rescaled Ricci harmonic map
heat flow.

There we prove short time existence using the same techniques as those
presented in \cite{ShiJDG1989,MilesC0,Riccistab}, see Theorem
\ref{exist thm}.  

In chapter \counterword{conv sec}, we show that the $L^2$-norm of
$g(t) -h$ is an exponentially decaying function of time (see Theorem
\ref{thm:intball}, Corollary \ref{intM cor}). This is the key
ingredient to the proofs of our stability results. The calculations to
prove this depend on an eigenvalue estimate for the Laplacian on
hyperbolic domains due to H. McKean \cite{McKean} and the closeness of
the evolving metric to that of hyperbolic space. In contrast to the
corresponding Euclidean result \cite{Riccistab}, we need strict
monotonicity of our integral quantity to establish long time
existence. The decay of the $L^2$-norm implies that the $C^0$-norm of
$g(t) - h$ is exponentially decaying in time (Theorem
\ref{thm:expdecayM}).  Interpolating between the $C^0$-norm and
$C^k$-norms, $k \in \N$, and using interior estimates, we see that all
of the $C^k$-norms are exponentially decaying in time (for $t \geq
1$).  This leads to long time existence and convergence.

In chapter \counterword{getting back sec}, we consider the related
scaled Ricci flow $\ti g(t)$ to the solution $g(t)$ obtained in
chapters \counterword{exist sec} and \counterword{conv sec}. The two
flows are related by time dependent diffeomorphisms $\phi_t:\H^n \to
\H^n$: $\ti g(t) :=\phi_t^* g(t)$.  As in the paper \cite{Riccistab},
we show that the estimates we obtained for $g(t)$ imply that $\ti g(t)
\to \psi^* h$ as $ t \to \infty$ in the $C^k$-norms. Here $\psi$ is a
diffeomorphism, and this diffeomorphism is the $C^k$-limit of the time
dependent diffeomorphisms $\phi_t$ which relate the two flows.  We
also show (as in \cite{Riccistab}) that $\psi \to \id$ as $|x| \to
\infty$, if the initial metric $g_0$ satisfies $g_0 -h \to 0$ as $x
\to \infty$ (see Theorem \ref{thm:backtoricci2}).  The proofs of this
chapter are the same (up to some minor modifications) as those of the
paper \cite{Riccistab}.

In Appendix \ref{scaling interior app} we gather various standard
results which we use in the paper.

In Appendix \ref{eucl space app} we show that the arguments used here
may be used in the Euclidean case to obtain analogous results (compare
with \cite{Riccistab}).

In Appendix \ref{conf sec} we address conformal Ricci flow in two
dimensions and obtain a result similar to the two-dimensional result
in \cite{Riccistab} without assuming that $|g-h|\to0$ near infinity. 

The authors were supported by the Deutsche Forschungsgemeinschaft,
DFG. 

\section{Existence}
\label{exist sec}
We first collect some evolution equations from \cite{ShiJDG1989}, and
then treat the question of short time existence.

In the following computations we always assume that in appropriate
coordinates, we have at a fixed point and at a fixed time
$h_{ij}=\delta_{ij}$ and $(g_{ij})=
\text{diag}(\lambda_1,\lambda_2,\ldots,\lambda_n)$, $\lambda_i >0$.

According to \cite[Lemma 2.1]{ShiJDG1989}, we get
\begin{align*}
  \delt g_{ij}=&\,g^{ab}\gradh_a\gradh_bg_{ij}
  -g^{kl}g_{ip}h^{pq}R_{jkql}(h)
  -g^{kl}g_{jp}h^{pq}R_{ikql}(h)\\
  &\,+\tfrac12g^{ab}g^{pq}\left(\gradh_i g_{pa}\gradh_jg_{qb}
    +2\gradh_ag_{jp}\gradh_qg_{ib}-2\gradh_ag_{jp}
    \gradh_bg_{iq}\right.\\
  &\,\left.\qquad\qquad\qquad-2\gradh_jg_{pa}\gradh_bg_{iq}
    -2\gradh_ig_{pa}\gradh_bg_{jq}\right)\\
  &\,-2(n-1)g_{ij}.\\
\end{align*}
Using that 
\[R_{ijkl}(h)=-(h_{ik}h_{jl}-h_{il}h_{jk})\] is the curvature tensor
of hyperbolic space of sectional curvature $-1$, we get
\begin{align*}
  -g^{kl}g_{ip}h^{pq}&R_{jkql}(h) -g^{kl}g_{jp}h^{pq}R_{ikql}(h)
  -2(n-1)g_{ij} \\
  = &\, g^{kl}g_{ip}h^{pq}(h_{jq}h_{kl} - h_{jl}h_{kq})
  +g^{kl}g_{jp}h^{pq}(h_{iq}h_{kl} - h_{il}h_{kq})
  -2(n-1)g_{ij}  \\
  = &\, 2\left(g^{kl}g_{ij}h_{kl} - h_{ij}\right) -2(n-1)g_{ij}\\
  = &\, 2\left(g^{kl}(h_{kl} -g_{kl})\right)g_{ij} +2(g_{ij} - h_{ij})\\
\end{align*}
and hence
\begin{lemma}\label{evolg lem}
  A metric $(g_{ij})$ solving \eqref{eq:RTFR} fulfills
  \begin{align*}
    \delt g_{ij}=&\,g^{ab}\gradh_a\gradh_bg_{ij}
    +  2g_{ij}\left(g^{kl}(h_{kl} -g_{kl})\right) +2(g_{ij} -
    h_{ij})\\ 
    &\,+\tfrac12g^{ab}g^{pq}\left(\gradh_i g_{pa}\gradh_jg_{qb}
      +2\gradh_ag_{jp}\gradh_qg_{ib}-2\gradh_ag_{jp}
      \gradh_bg_{iq}\right.\\
    &\,\left.\qquad\qquad\qquad-2\gradh_jg_{pa}\gradh_bg_{iq}
      -2\gradh_ig_{pa}\gradh_bg_{jq}\right).\\
  \end{align*}
\end{lemma}

For later use, we derive the evolution equation of $|g-h|^2$.
\begin{lemma}\label{evolZ lem}
  Let $g\in\M\infty{\H^n,(0,T)}$, $T>0$, be a solution to
  \eqref{eq:RTFR} which is $\epsilon$-close to the hyperbolic metric
  $h$ of sectional curvature $-1$. Assume that $\epsilon>0$ is
  sufficiently small. Then
  \begin{equation}
    \delt |g-h|^2\le
    g^{ij}\nabla_i\nabla_j|g-h|^2
    -(2-\epsilon)\left|\nabla(g-h)\right|^2 +(4+\epsilon)|g-h|^2, 
  \end{equation}
  where we write $\epsilon$ instead of $c(n)\epsilon$ and $\nabla$
  instead of $\gradh$. 
\end{lemma}
\begin{proof}
  Note that the norm of a tensor $Z$ of order $m$ fulfills
  $$\uph|Z|^2\equiv|Z|^2\le(1+\epsilon)\cdot\up g|Z|^2
  \le(1+\epsilon)\cdot\uph|Z|^2.$$ Choose coordinates such that
  $h_{ij}=\delta_{ij}$ and $g_{ij}=\lambda_i\delta_{ij}$. We use $*$
  similarly as in \cite[Ch.{} 13]{HamiltonThree} to denote
  contractions with respect to $h$, $g$ or their inverses. Let
  $Z=g-h$. Lemma \ref{evolg lem} yields
  \begin{align*}
    \partt |Z|^2 \equiv&\, \partt |g-h|^2 = 2\sum\limits_i(g_{ii} -
    h_{ii})\left(\partt g_{ii}\right)\\
    =&\, g^{ij} \grad_i \grad_j |g-h|^2 - (2-\epsilon) |\grad g|^2\\
    &\,+ 2\sum\limits_i (g_{ii} - h_{ii})\left[2(g_{ii} - h_{ii})
      -2g_{ii}\sum\limits_k\left(g^{kk}(g_{kk} -h_{kk})\right)
    \right]\\
    &\,+ \sum\limits_i(g_{ii} - h_{ii}) (\grad g * \grad g)_{ii} \\
    \leq&\,  g^{ij} \grad_i \grad_j |g-h|^2 - (2- \ep) |\grad (g -h) |^2\\
    &\, + 4\sum\limits_i(g_{ii} - h_{ii})\left[(g_{ii} - h_{ii})
      -g_{ii}\sum\limits_k\left(g^{kk}(g_{kk} -h_{kk})\right) \right].
  \end{align*}
  Let us examine the zeroth order term $S$ of the above equation.
  \begin{align*}
    S=&\, 4\sum (g_{ii} - h_{ii})\left[(g_{ii} - h_{ii})
      -g_{ii}\sum\left(g^{kk}(g_{kk}
        -h_{kk})\right) \right]\\
    =&\, 4 \sum_i (\lambda_i-1)^2-4 \sum_i \lambda_i(\lambda_i-1)
    \sum_k\left(1-\frac1{\lambda_k}\right) \\
    \leq&\, (4+\ep) |Z|^2-4 \left(\sum_i (\lambda_i-1)\right)^2 .
  \end{align*}
  The claim follows. 
\end{proof} 

We use this to show that we can solve Dirichlet problems for
\eqref{eq:RTFR} on a short time interval. In the following we pick a
point $p_0$ and fix it throughout. For simplicity of notation we will
denote this point with $0$. All balls $B_R(0)$ are geodesic balls with
respect to the hyperbolic metric $h$.

The following result also extends readily to \eqref{eq:RTFR} on all of
$\H^n$ provided that a non-compact maximum principle is applicable.
\begin{corollary}
  \label{close cor}
  Let $g\in\mathcal M^\infty_0(B_R(0),[0,T)), \ 0<T<\infty,$ be a solution to
  \eqref{eq:RTFR} on $B_R(0)\subset\H^n$ instead of $\H^n$ with
  $g(t)|_{\partial B_R(0)}=h|_{\partial B_R(0)}$. Let $0<\delta$. Then
  there exists $\epsilon=\epsilon(n,T,\delta)>0$ such that
  $\sup_{B_R(0)}|g(0)-h|\le\epsilon$ implies 
  \[\sup\limits_{B_R(0)\times[0,T)}|g-h|\le\delta.\]
\end{corollary}
\begin{proof}
  Assume without loss of generality that $\delta$ is smaller than
  $\epsilon$ in Lemma \ref{evolZ lem} and $\epsilon\le1$. Fix
  $\epsilon:=\delta e^{-5T}$. Then Lemma \ref{evolZ lem} implies that
  $\sup_{B_R(0)}|g(t)-h|\le\delta e^{-5(T-t)}$ as long as
  $\sup_{B_R(0)}|g(t)-h|\le\epsilon$. The result follows. 
\end{proof}

If solutions to \eqref{eq:RTFR} stay sufficiently close to the
hyperbolic metric $h$ of sectional curvature $-1$, they exist for all
times.
\begin{theorem}
  \label{exist thm}
  For all $n\in\N$ there exists a $\ti \de(n) >0$ such that the
  following holds. Let $0<\epsilon<\delta\le\ti \de (n)$.  Then every
  metric $g_0\in\M0{\H^n}$ with $\Vert
  g_0-h\Vert_{L^\infty}\le\epsilon$ has a $\de$-{\rm maximal} solution
  $g\in\mathbb M^\infty_0(\H^n,[0,T_{g_0}))$ to \eqref{eq:RTFR}, where
  $T_{g_0} >0$ and $\Vert g(t)-h\Vert_{L^\infty}<\delta$ for all $t
  \in [0,T_{g_0})$. The solution is $\de$-{\rm maximal} in the
  following sense.  Either $T_{g_0}=\infty$ and $\Vert
  g(t)-h\Vert_{L^\infty}<\delta$ for all $0\le t$ or we can extend $g$
  to a solution on $\H^n\times[0,T_{g_0}+\tau)$, for some $\tau=
  \tau(n)>0$, and $\Vert g(T_{g_0})-h\Vert_{L^\infty}=\delta$.
\end{theorem}
\begin{proof}
  The proof follows directly from the corresponding proofs in
  \cite{MilesC0,ShiJDG1989,Riccistab}: We mollify $g_0$ and obtain
  $g^i_0$, $i\in\N$.  Then we consider $g^{i,R}_0:=\eta
  g^i_0+(1-\eta)h$, $R\ge1$, where $\eta:\H^n\to\H^n$ is a smooth
  function fulfilling $\eta=1$ in $B_R(0)$, $\eta=0$ outside
  $B_{2R}(0)$ and $|\nabla\eta|\le2/R$.  Hence $\big\Vert
  g^{i,R}_0-h\big\Vert_{L^\infty} \le\left\Vert
    g^i_0-h\right\Vert_{L^\infty} \leq \Vert
  g_0-h\Vert_{L^\infty(\H^n)} \leq \epsilon $.  Arguing as in
  \cite{MilesC0,ShiJDG1989,Riccistab} (see Section 3 of
  \cite{ShiJDG1989}), and using that $g^{i,R}_0$ is $\ti \de(n)$ close
  to $h$, we see that there exist solutions
  $g^{i,R}\in\M\infty{B_{3R}(0),[0,\tau]}$ with
  $0<\tau=\tau(n)<\infty$ to \eqref{eq:RTFR} on $B_{3R}(0)$ with
  $g^{i,R}(0)=g^{i,R}_0$ on $B_{3R}(0)$ and $g^{i,R}=h$ on $\partial
  B_{3R}(0)\times[0,\tau]$.  From Lemma 5.1 of \cite{MilesC0} we see
  that we have interior estimates of the form $|\grad^j g^{i,R}|^2
  \leq c/t^j$ on balls of radius $R/2$ for all $t \in [0,\tau]$ for
  all $j \in \N$.  Taking a diagonal subsequence and using the Theorem
  of Arzel\`a-Ascoli, we obtain a solution $g \in
  \M\infty{\H^n,(0,\tau]}$. Furthermore, $g(t) \to g_0$ as $t \to 0$
  in the $C^0_{\loc}$-norm as we see from Theorem 5.2 in
  \cite{MilesC0}.

  If $\Vert g(t) - h\Vert_{L^\infty} < \de$ for all $t \in [0,\tau]$,
  then we may repeat this argument in view of the fact that $\de\le \ti
  \de(n)$.  By induction we obtain a solution $g \in
  M^\infty_0(\H^n,[0,S))$ where either
  \begin{enumerate}[(a)]
  \item\label{S infty} $S = \infty$ and $\Vert g(t) -
    h\Vert_{L^\infty} < \de$ for all $t>0$, or
  \item\label{S finite} $S =N \tau$ for some $N \in \N$ and $\Vert g(t) -
    h\Vert_{L^\infty} < \de$ for all $t \in [0,(N-1)\tau]$ but there
    exists at least one time $\ti t\in ((N-1)\tau,N\tau]$ with $\left\Vert
      g\left(\ti t\right) - h\right\Vert_{L^\infty} = \de$.
  \end{enumerate}
  In the case \eqref{S infty}, we set $T_{g_0} = \infty$ and we are
  finished.

  So assume we are in case \eqref{S finite} for the rest of the
  argument.  In view of Lemma \ref{evolZ lem}, the maximum principle,
  and the way we constructed our solutions, we see that in fact there
  is a first time $T_{g_0} \in ((N-1)\tau,N\tau]$ with $\Vert
  g(T_{g_0}) - h\Vert_{L^\infty} = \de$ and $\Vert g(t) -
  h\Vert_{L^\infty} < \de$ for all $t \in [0,T_{g_0})$.  Using that
  $g(T_{g_0})$ is $\de$ close to $h$ and $\de \leq \ti \de$ we may
  repeat the first part of the proof to obtain a solution defined on
  $[0,T_{g_0} + \tau]$.  This completes the proof.
\end{proof}

\begin{proposition}
  \label{lower exist time}
  Assume everything is as in Theorem \ref{exist thm}, and $\hat T>0$ be
  given.  If we choose $\ep = \ep\big(n,\delta,\hat T\big)>0$
  small enough in the above Theorem, then the solution $g\in\mathbb
  M^\infty_0(\H^n,[0,T_{g_0}+\tau))$ from Theorem \ref{exist thm}
  satisfies $T_{g_0} \geq \hat T$.
\end{proposition}
\begin{proof}
  By Corollary \ref{close cor}, we can choose
  $\ep=\epsilon\big(n,\delta,\hat T\big)$ small enough so that $\Vert
  g(t)-h\Vert_{L^\infty} < \de$ for all $t \in \big[0,\hat
  T\big]\cap[0,T_{g_0}]$. Theorem \ref{exist thm} yields a
  contradiction if $T_{g_0}<\hat T$.
\end{proof}

\section{Convergence}
\label{conv sec}
Convergence is based on a Lyapunov function.
\begin{theorem}\label{thm:intball} 
  Let $n\geq 4$.  There exists $\delta_0=\delta_0(n)>0$ such that the
  following holds. Let $g\in\mathcal{M}^\infty(B_{R}, [0,T))$ be a
  solution to \eqref{eq:RTFR} with $g=h$ on $\partial B_R(0)\times
  [0,T)$ and assume that $\sup_{B_R(0)\times[0,T)}|g-h|\leq
  \delta_0$. Then we have
  \[\int\limits_{B_R(0)}|g(t)-h|^2 \, d\!\vol_h \leq e^{-\alpha t}
  \int\limits_{B_R(0)}|g(0)-h|^2 \, d\!\vol_h\] 
  for $\alpha(n):=(2(n-1)^2-17)/4\geq \tfrac{1}{4}$.
\end{theorem}

\begin{proof} 
  Assume that $\delta_0$ is such that $g$ is
  $\epsilon=\epsilon(n)$-close to $h$ for some sufficiently small
  $\epsilon>0$.  We compute, using Lemma \ref{evolZ lem},
  \begin{align*}
    \partt \int\limits_{B_R(0)} |Z|^2\, d\!\vol_h \leq&\,
    \int\limits_{B_R(0)} g^{ij} \grad_i \grad_j |Z|^2 - (2-\ep)|\grad
    Z|^2 + (4 + \ep)|Z|^2\,d\!\vol_h\\
    =&\, \int\limits_{\partial B_R(0)} \nu_i g^{ij} \grad_j |Z|^2 
    -\int\limits_{B_R(0)} \grad_jg^{jk} \grad_k|Z|^2\,d\!\vol_h\\
    &\, +\int\limits_{B_R(0)}- (2-\ep)|\grad Z|^2 + (4 +
    \ep)|Z|^2\,d\!\vol_h\\
    \leq&\, \int\limits_{B_r(0)} - (2-\ep)|\grad |Z||^2 + (4 +
    \ep)|Z|^2\,d\!\vol_h\ ,
  \end{align*}
  where we used that $\left|\grad_ig^{ij} \grad_j|Z|^2\right|\leq \ep
  |\nabla Z|^2$ and that on $\partial B_R(0)$ the gradient $\nabla
  |Z|^2$ is anti-parallel to the outer unit normal $\nu$.  Furthermore
  we apply Kato's inequality $ |\grad |Z||^2 \leq |\grad Z|^2$ which
  is valid whenever $|Z| \neq 0$ and for Sobolev functions.

  Using McKean's inequality \cite{McKean} for the first eigenvalue
  $$\sigma_1 \geq \frac{(n-1)^2}{4}$$ 
  on hyperbolic domains we see
  \begin{align*}
    \partt \int\limits_{B_R(0)} |Z|^2\, d\!\vol_h \leq&\,
    \int\limits_{B_R(0)} - (2-\ep)|\grad |Z||^2 + (4 +
    \ep)|Z|^2\,d\!\vol_h\\
    \leq&\, \frac{8 - (n-1)^2+ \ep}{2}\ \int\limits_{B_R(0)}|Z|^2\,
    d\!\vol_h\ .
  \end{align*}
  Assuming that $\ep < 1/2$, we can choose
  $$ \alpha:= (2(n-1)^2-17)/4\ .$$
\end{proof}

Since for the proof of existence of a solution to \eqref{eq:RTFR} we
use Dirichlet problems as above, this monotonicity extends to the
constructed solutions on $\H^n\times [0,T)$: Let $g^{i,R}(t)$ be as in
Theorem \ref{exist thm}. Then we get 
\[\left\Vert g^{i,R}(t)-h\right\Vert^2_{L^2(B_{3R}(0))}
\le e^{-\alpha t}\big\Vert g^{i,R}_0-h\big\Vert^2_{L^2(B_{3R}(0))}
\le e^{-\alpha t}\Vert g_0-h\Vert^2_{L^2(\H^n)}\]
As $R\to\infty$, we obtain
\begin{corollary} \label{intM cor} Let $n\ge4$ and $T>0$ be
  given. Assume that $g_0\in\M\infty{\H^n}$ satisfies $\Vert
  g_0-h\Vert_{L^2(\H^n)}<\infty$. Then there exists
  $\epsilon_0=\epsilon_0(n,T)$ such that, if $\sup_{\H^n} |g_0-h|\leq
  \epsilon_0$ then a solution $g\in\mathcal{M}^\infty(\H^n, [0,T))$ to
  \eqref{eq:RTFR} with $g(\cdot,0)=g_0(\cdot)$ exists and
  $\sup_{\H^n\times[0,T)}|g-h| \leq \delta_0$, where $\delta_0$ is as
  in Theorem \ref{thm:intball}. Furthermore we have the estimate
  \[\Vert g(t)-h\Vert^2_{L^2(\H^n)} \le e^{-\alpha t}
  \Vert g_0-h\Vert^2_{L^2(\H^n)}\]
  for all $t \in [0,T)$, where $\alpha = \alpha(n)\geq\tfrac{1}{4}$.
\end{corollary}
\begin{proof}
  Existence and closeness to $h$ follow from Corollary \ref{close
    cor}, Proposition \ref{lower exist time} and Theorem \ref{exist thm}.
\end{proof}

Using the gradient estimate we see that the exponential convergence of
the $L^2$-norm of $|g-h|$ also implies exponential convergence in the
$\sup$-norm, compare \cite[Lemma 7.1]{Riccistab}.

\begin{theorem}\label{thm:expdecayM}
  Let $n\ge4$.  Assume that $g \in \mathcal{M}^\infty(\H^n, [0,T))$ is
  a solution to \eqref{eq:RTFR} with $\Vert g(0)-h\Vert_{L^2(\H^n)} =:
  K < \infty$, $\sup_{{\H^n}\times[0,T)}|g-h| \leq \delta_0$ and
  \[\Vert g(t)-h\Vert^2_{L^2(\H^n)}\leq e^{-\alpha t}
  \Vert g(0)-h\Vert^2_{L^2(\H^n)}\ , \] where $\delta_0$ is as in
  Theorem \ref{thm:intball}. Then
  \begin{equation} \label{exp} 
    \sup_{{\H^n}} |g(t)-h|\leq C(n,K) e^{-\beta t},
  \end{equation}
  where $\beta=\frac\alpha{n+2} =\frac{2(n-1)^2-17}{4(n+2)}>0$.
\end{theorem}

\begin{proof} 
  We can assume w.l.o.g.{} that $\delta_0< 1$. We choose $\tau:=
  \tfrac{n+1}{\alpha}\ln(\delta_0^{-1})>0$. Note that this implies
  $$ \sup_{{\H^n}}|g(t)-h|\leq e^{-\beta t}$$
  for $t\in[0,\tau)$ and $\beta:= \alpha/(n+1)$. By the interior
  estimates of the form $\big|\gradh g(t)\big|\le c\cdot t^{-1/2}$,
  there exists a constant $C^\prime=C^\prime(n)$, such that
  $$ \left|\gradh g(\cdot,t)\right|^h \leq C^\prime $$
  for $t\in [\tau, T)$. Fix such a $t \in [\tau,T)$. Let $\gamma:=
  \sup_{\H^n}|g(t)-h|$ and choose a point $p_0\in {\H^n}$ such that
  $|g(p_0,t)-h(p_0) |\geq \tfrac{1}{2} \gamma$. By the gradient
  estimate, we have
  $$ |g(\cdot,t)-h|\geq \frac{1}{4} \gamma $$
  on $B_{\gamma/(4C^\prime)}(p_0)$. This implies
  $$\Vert g(t)-h\Vert^2_{L^2(\H^n)} \geq \omega_n
  (C^\prime)^{-n}\Big(\frac{\gamma}{4}\Big)^{n+2}\ ,$$ where
  $\omega_n$ is the measure of the unit ball in $\R^n$. This yields
  $$ \gamma \leq 4 (C^\prime)^\frac{n}{n+2}
  \left(\frac{K^2}{\omega_n}\right)^{\frac{1}{n+2}}
  e^{-\frac{\alpha}{n+2} t}\ .$$ Choosing $C(n,K) = 1+ 4
  (C^\prime)^\frac{n}{n+2}
  \big(\frac{K^2}{\omega_n}\big)^{\frac{1}{n+2}}$ we have \eqref{exp}.
\end{proof}

This $\sup$-estimate allows us to construct a solution which exists
for all times.

\begin{theorem}\label{thm:existenceBIG}
  Let $n\ge4$. For all $K>0$ there exists
  $\epsilon_1=\epsilon_1(n,K)>0$ such that the following holds. Let
  $g_0\in\M\infty{\H^n}$ satisfy $\Vert g_0-h\Vert_{L^2(\H^n)}\leq K$
  and $\sup_{\H^n}|g_0-h| \leq \epsilon_1$. Then there exists a
  solution $g \in \mathcal{M}^\infty({\H^n}, [0,\infty))$ to
  \eqref{eq:RTFR} with $g(0)=g_0$ such that
  \begin{equation}
    \label{gh exp decay eq}
    \sup_{\H^n} |g(t)-h| \leq C(n, K) e^{-\beta t}
  \end{equation}
  for $\beta=\beta(n)$ as in Theorem \ref{thm:expdecayM}.
\end{theorem}

\begin{proof}
  According to Theorem \ref{exist thm}, we obtain existence for all
  times if we can prove the estimate $\Vert g(t)-h\Vert_{L^\infty}\le
  \ti \delta= \ti \delta_{\text{Thm.{} \ref{exist thm}}}$ for all $t$
  for any a priori solution (that is, we must prove the estimate for
  all $t$ that the solution is defined).  Given any $T>0$, we can
  choose $\ep(n,T)>0$ small enough so that such an estimate will hold,
  in view of Proposition \ref{lower exist time} and Theorem \ref{exist
    thm} for $0\le t<T$.  Theorem \ref{thm:intball} implies integral
  bounds which combined with Theorem \ref{thm:expdecayM} yields
  $|g(t)-h|\le\tilde\delta$ if $t\ge T$ and $T$ is chosen sufficiently
  large.  Choose $T$ and $\ep(n,T)$ so that both of these requirements
  are satisfied.

  This implies long time existence.\par Theorem
  \ref{thm:expdecayM} also implies \eqref{gh exp decay eq} for $t\ge
  T$. Fixing $C(n,K)$ such that $C(n,K)\le\delta\cdot e^{\beta T}$
  we obtain \eqref{gh exp decay eq} for all times.
\end{proof}

By interpolation the exponential decay extends to higher derivatives
of the evolving metric.

\begin{theorem}\label{thm:gradest}
  Let $n\ge4$.  Let $g_0\in\M\infty{\H^n}$ and $g \in
  \mathcal{M}^\infty({\H^n}, [0,\infty))$ be as in Theorem
  \ref{thm:existenceBIG}. We have additionally
  $$ \sup_{\H^n} \Big|\gradh^j g(t)\Big| \leq C(n,j,K,(\beta_j))
  e^{-\beta_j t}$$ where $0<\beta_j<\beta(n)$, $\beta(n)$ as in
  Theorem \ref{thm:expdecayM}, is arbitrary.  In particular,
  \[\lim_{t \to \infty} \sup_{\H^n} \Vert g(t) -h\Vert_{C^k({\H^n})} =
  0,\] where $\Vert S\Vert_{C^k} := \sum\limits_{i = 0}^n
  \sup\limits_{\H^n} |\grad^j S|^2$.
\end{theorem}
\begin{proof}
  From the interior estimates in Lemma \ref{lem:intest}, we have
  $\sup_{{\H^n}} |\grad^j g|^2(t) \leq {c(n,j)/ (t-L)},$ for all $t
  \in [L,L+1]$. In particular, $\sup_{{\H^n}} |\grad^j g|^2(L+1) \leq
  c(n,j).$ Hence, as $L>0$ was arbitrary, we get
  \begin{eqnarray}\label{gradint}
    \sup_{{\H^n}} \left|\grad^j g\right|^2(\cdot,t) \leq c(n,j)
  \end{eqnarray}
  for all $t\geq 1$.  Interpolating on a ball of radius one (see Lemma
  \ref{lemm:interp}) gives us
  $$\sup_{\H^n} \left|\grad^j g\right|^2(t) \leq  \ti c(n,j)
  \left(\sup_{\H^n} |g(t)-h|\right)^{\frac1{2^{j-1}}} \leq C(n,j,K)
  e^{-\frac{\beta}{2^{j-1}} t}$$ in view of \eqref{gradint} and
  \eqref{exp}. Iterated interpolation yields the result, see e.\,g.{}
  \cite{Riccistab}.
\end{proof}

\begin{proof}[Proof of Theorem \ref{thm:main}]
  As the decay of $|g(t)-h|$ as $t\to\infty$ obtained in this section
  does not depend on the smoothness of $g_0$, we can approximate
  $g_0\in\M0{\H^n}$ and pass to a limit to obtain Theorem \ref{thm:main}.
\end{proof}

\begin{proof}[Proof of Theorem \ref{thm:sol given}]
  Local closeness estimates (see Lemma \ref{lem:intcloseness})
  show that 
  \[\lim\limits_{r\to\infty}\Vert g(t)-h\Vert_{L^\infty(\H^n\setminus
    B_r(0))}=0\] is preserved during the flow, even uniformly on
  compact time intervals. Hence 
  \[\max\left\{|g(t)-h|^2-\delta,0\right\}
  \equiv\left(|g-h|^2-\delta\right)_+\] has compact support on $\H^n
  \times[0,K]$ for all $K< \infty$ and we may consider the integral
  $I_\delta:=\big\Vert\left(|g-h|^2-\delta\right)_+\big
  \Vert_{L^1\left(\H^n\right)}$ for any $\delta>0$, which is similar
  to the integral $I^{m,p}_\delta$ defined in \cite[Theorem
  6.1]{Riccistab} or to $I^p_\delta$ in Appendix \ref{eucl space
    app}. The techniques of the proof of Theorem \ref{thm:intball} and
  approximations as in \cite[Theorem 6.1]{Riccistab} imply for $R\gg1$
  that
  \[I^R_\delta(t):=\big\Vert
  \left(|g(t)-h|^2-\delta\right)_+\big\Vert_{L^1(B_R(0))}\le
  e^{-\alpha t}\cdot I^R_\delta(0)\le e^{-\alpha t}\cdot\Vert
  g(0)-h\Vert^2_{L^2(\H^n)}. \] The rest of the proof is similar to
  the proof of Theorem \ref{thm:main}.
\end{proof}

\section{Getting back to Ricci Flow}
\label{getting back sec}
\begin{theorem}\label{thm:backtoricci}
  Assume that $g \in \mathcal{M}^\infty({\H^n}, [0,\infty))$ is the
  solution to \eqref{eq:RTFR} coming from Theorem
  \ref{thm:existenceBIG}.  Then there exists a smooth map $\phi:{\H^n}
  \times [0, \infty) \to {\H^n} $ such that $\phi(\cdot,t) =: \phi_t:
  {\H^n} \to {\H^n}$ is a diffeomorphism, $\phi_0 = \id$ and $\tilde
  g(\cdot,t):=(\phi_t)^*g(\cdot,t)$ is a smooth solution to the {\it
    scaled Ricci flow}
  \begin{equation*}
    \partt g = -2 \Ricci - 2(n-1)g 
  \end{equation*}
  with $\tilde g_t \to g_0$ as $t \searrow 0$.  Furthermore there
  exists a smooth diffeomorphism $\psi: {\H^n} \to {\H^n}$ with
  $\phi_t \to \psi$ as $t \to \infty$ and $\tilde g_t \to \psi^* g_0$
  as $t \to \infty$.  Here convergence is in $C^k$ on $\H^n$ for all $k$.
\end{theorem}
\begin{proof}
  This argument is the same as in Lemma 9.1 and Theorem 9.2
  of \cite{Riccistab} with some
  minor differences.  We explain here where the argument of
  \cite{Riccistab} must be modified in order for it to work in this
  case.

  As explained in Lemma 9.1 in \cite{Riccistab}, we can construct 
  smooth maps $\phi:\H^n
  \times [0,\infty) \to \H^n$ such that
  $$\begin{cases}
    \delt {\phi}^{\al}(x,t) ={V}^{\al}(\phi(x,t),t),&
    (x,t)\in \H^n \times [0,\infty),\\
    {\phi}(x,0) = x,&x\in \H^n,
  \end{cases}$$ where $ {V}^{\al}(y,t) := -{g}^{\be \ga} \left({\up{g}
      \Gamma}^{\al}_{\be \ga} - {\uph \Gamma}^{\al}_{\be \ga}\right)
  (y,t)$ and $\phi_t:= \phi(\cdot,t):\H^n \to \H^n$ are
  diffeomorphisms. (Compared to \cite{Riccistab}, we have changed the
  sign in the definition of $V$ in order to correct a typo there.)
  This part of the proof is the same. A direct calculation shows that
  $\ti g(t):= \phi_t^*g(t)$ solves the scaled Ricci flow equation.

  In Theorem 9.2 of \cite{Riccistab} it is shown that $\phi_t \to
  \phi_{\infty}$ as $t \to \infty$ where $\phi_{\infty}:\H^n \to \H^n$
  is a smooth diffeomorphism, and the convergence is in $C^k$ (for all
  $k$) on $\H^n$.  The proof of this is carried out in three steps.

  In step 1 it is shown that $|\partt \phi_t(x)| \leq \frac{1}{
    t^{r}}$ for some $r>1$, for all $t \geq 1$, and $|\phi_t(x) - x|
  \leq c$ for all $t$.

  In step 2, the existence of a smooth function $\phi_{\infty}: \H^n
  \to \H^n$ with $\phi_t \to \phi_{\infty}\equiv\psi$ as $t \to
  \infty$ is shown.

  The proofs of steps one and two carry over to this situation
  without any changes.

  In step 3, it is shown that $\phi_{\infty}$ is a diffeomorphism.
  This proof carries over with some minor modifications which we
  describe in the rest of the proof here.

  Letting $l(t):= \phi_t^{*} g(t),$ we know that $l$ solves the {\it
    scaled Ricci flow} \eqref{eq:RFresc} on $\H^n$, and that
  $$\sup_{\H^n} {}^{l(t)}|\Ricci(l(t)) +  2(n-1)l(t)| = 
  \sup_{\H^n} {}^{g(t)}|\Ricci(g(t)) + 2(n-1)g(t)| \leq e^{-\be t}.$$
  for all $t>0$ for a $\be >0$, in view of Theorem \ref{thm:gradest}.
  Hence $l(t)$ converges locally uniformly (smoothly) to a smooth
  metric $l_{\infty}$ on $\H^n$ as explained in \cite{Riccistab}.

  Choose geodesic coordinates for $h$ centred at $y$ in $B_{1}(y)$.
  Now using the definition of ${l}$, and the uniform convergence of
  $l$ we get
  \begin{align*}
    {\frac 1 c} \de_{\al \be} \leq l_{\al \be}(x,t)
    &\,= \parti{\phi^s_t}{y^{\al}}(x,t)
    \parti{\phi^k_t}{y^{\be}}(x,t)
    \,{g}_{s k}(\phi_t(x),t) \\
    &\,\leq (1 +\tilde\ep) \parti{\phi^s_t}{y^{\al}}
    \parti{\phi^k_t}{y^{\be}}(x,t) h_{sk}(\phi_t(x))\\
    &\,\leq c (1 +\tilde\ep)\left(D\phi_t\right)
    \left(D\phi_t\right)^T(x,t).
  \end{align*}
  In particular, we see that $\det\left(D\phi_t\right)^2(x)
  \geq\frac1{(1 +\tilde\ep)^nc}>0$ for $x \in B_1(y)$ , where $Df$ is
  the Jacobian of $f$.  As explained in \cite{Riccistab}, this shows
  that $\phi_\infty$ is a diffeomorphism.
\end{proof}

\begin{theorem}\label{thm:backtoricci2}
  Let everything be as in the above Theorem \ref{thm:backtoricci},
  with the extra assumption that $\sup_{({\H^n}\setminus B_r(0))}
  |g_0-h| \to 0$ as $r \to \infty$.  Then the diffeomorphism $\psi$
  appearing in the above Theorem satisfies $\psi \to \id$ as $x \to
  \infty$ (in $C^k$ for all $k$).  In particular, for every
  $\eta>0$, there exists an $R>0$ such that
  $$\sup_{\H^n \backslash B_R(0)} |\phi_t(x) - x | \leq \eta$$
  for all $t$.
\end{theorem}
\begin{proof} 
  The proof is completely analogous to the proof of Lemma 9.3 in
  \cite{Riccistab}. Let $\eta >0$.  From Lemma \ref{lem:intcloseness}
  and the estimates of Theorem \ref{thm:existenceBIG} we can choose
  $R>0$ large so that
  $$ | g(t) -h| \leq \eta \mbox{ on } \H^n \setminus B_R(0)$$ for all
  $t>0$. 
  From the interior estimates of \cite{MilesC0} (see Lemma
  \ref{lem:intest}) we get
  $$  |\grad ^2 g| \leq  \frac c  t$$ for $t \in [0,1]$, and
  hence
  $$  |\grad ^2 g| \leq  \frac c  t$$ for all $t \in [0, \infty)$,
  in view of Theorem \ref{thm:existenceBIG} and interpolation with
  higher order derivatives, see Lemma \ref{lem:intest} and Lemma
  \ref{lemm:interp}.  Interpolating between the $C^0$-norm and the
  $C^2$-norm on a ball of radius one for $t \leq 1$ (see Lemma
  \ref{lemm:interp}) we get
  $$  |\grad g| \leq  
  \frac {\sqrt c \sqrt \eta} {\sqrt t}$$ for all $t \in [0, \infty)$
  on $\H^n \setminus B_R(0)$.  Arguing as in Step 1 of the above
  Theorem gives us the result.
\end{proof}

\begin{appendix}
  \section{Scaling and Interior Estimates}
  \label{scaling interior app}

  \begin{lemma}\label{lem:intest}
    Let $(\H^n,g(t))_{t\in[0,T)}$ be a solution to equation
    \eqref{eq:RTFR}, with $|g(t)-h|\le\epsilon(n)$ for $\epsilon(n)>0$
    small enough. Then \[\sup\limits_{\H^n}
    \Big|\gradh^jg(\cdot,t)\Big|^2 \le \frac{c(j,n)}{t^j}\] for all
    $t\le\min\{1,T\}$. 
  \end{lemma}
  \begin{proof}
    This is Theorem 4.3 in \cite{MilesC0}.
  \end{proof}

  \begin{lemma}\label{lem:intcloseness} 
    Fix a point $p_0 \in {\H^n}$ and let $g \in
    \mathcal{M}^\infty({\H^n},[0, T))$, where $T \in (0,\infty]$,
    be a solution to
    \eqref{eq:RTFR} such that
    $$\sup_{{\H^n}\setminus B_r(p_0)}|g(\cdot,0) -h|
    \rightarrow 0 \ \text{as}\ r\rightarrow \infty. $$ Then for every
    $0<\tau<T$ and $0<\ep\leq 1$ there exist an $R_0>0$ such that
    $$ \sup_{({\H^n}\setminus B_{R_0}(p_0))\times [0,\tau]}|g-h| \leq
    \ep\ .$$  
  \end{lemma}
  \begin{proof} 
    Choose a smooth function $\eta:\R \rightarrow \R, 0\leq \eta\leq
    1$, such that $\eta \equiv 1$ on $(-\infty,1]$, $\eta \equiv 0$ on
    $[2,\infty)$ and $\eta^\prime \leq 0$. We can furthermore assume
    that $|\eta^{\prime\prime}|\leq 8$ and $|\eta^\prime|^2 \leq 16
    \eta$. \\
    Let $\rho_{p}(\cdot)$ denote the distance to a point $p \in
    {\H^n}$ with respect to the hyperbolic metric, and define the
    cut-off function
    $$ \gamma_{p,R}:= \eta\Big(\frac{\rho_{p}}{R}\Big)\ .$$
    Then we have, suppressing in the following the subscripts $p$ and
    $R$,
    $$|\nabla\gamma| \leq \frac{|\eta^\prime|}{R}, \qquad |\nabla^2\gamma|
    \leq \frac{C}{R}+\frac{C}{R^2}\ .$$ Define $\psi:=\gamma
    |Z|^2$. Using Lemma \ref{evolZ lem} and the above estimates we see
    that
    \begin{align*}
      \partt \psi \leq&\, g^{ij}\nabla_i\nabla_j \psi - 2
      g^{ij}\nabla_i \gamma \nabla_j |Z|^2 - |Z|^2
      g^{ij}\nabla_i\nabla_j \gamma\\
      &\, -(2-\ep)\gamma |\nabla Z|^2 + (4+\ep) \psi\\
      \leq&\, g^{ij}\nabla_i\nabla_j \psi +(4+\ep)\psi +\frac{\ep}{R}+
      \frac{\ep}{R^2}\ ,
    \end{align*}
    where we used the estimates on the derivatives of $\gamma$ and
    Kato's inequality. Note that the closeness assumption can be
    justified inductively using continuity on small time intervals and
    the following argument. Now choosing $\text{dist}(p,p_0)$ big
    enough and $R=\tfrac{1}{2}\text{dist}(p,p_0)$ we can integrate the
    inequality to yield the desired estimate.
  \end{proof}

  \begin{lemma}
    The scaled Ricci flow and the scaled Ricci harmonic map heat flow
    are related as follows.  Assume $\phi_t:{\H^n} \to {\H^n}$ solves
    \[\delt\phi_t(x,t)=-V(\phi_t(x,t),t)\,\] where the components of
    $V$ are given by $V^\alpha:=g^{\beta\gamma}\left(\up
      g\Gamma^\alpha_{\beta\gamma}-\up
      h\Gamma^\alpha_{\beta\gamma}\right)$ and that the $\phi_t:{\H^n}
    \to {\H^n}$ are smooth and diffeomorphisms for all time.  Let
    $({\H^n},g(t))_{t \in [0,T)}$ be a solution to the scaled Ricci
    harmonic map heat flow \eqref{eq:RTFR}, $V_i=g_{i\alpha}V^\alpha$.
    Then $({\H^n},\ti g(t))_{t\in [0,T)}$ solves the scaled Ricci flow
    \eqref{eq:RFresc}, with $\ti g(0) = g(0)$, where here $\ti g(t):=
    \phi_t^*(g(t)).$
  \end{lemma}
  \begin{proof}
    For $\ti g(t):= \phi_t^* g(t)$, we get
    \begin{align*}
      \partt (\ti g(t)) = &\,(\phi_t)^* \left(\partt g\right) +
      \left.\parts\right|_{s = 0}(\phi^*_{t+s} g(t)) \\
      = &\,-2\Ricci(\ti g(t)) -2(n-1) \ti g(t) +
      \phi_t^*(\curlL_{V(t)}
      g(t)) - \curlL_{ (\phi_t^{-1})_{*} V(t)} (\phi_t^*
      g(t))\ \\ 
      = &\,-2\Ricci(\ti g(t)) -2(n-1) \ti g(t),
    \end{align*}
    where here $\curlL_W k $ is the Lie-derivative of $k$ in the
    direction $W$ (in coordinates $(\curlL_W k)_{ij} = \up{k} \grad_i
    W_j + \up{k} \grad_j W_i$), see \cite[Chapter 2, Section
    6]{ChLuNi}.
  \end{proof}

  \begin{lemma}\label{lem:equiv}
    The Ricci flow
    \[\partt g = -2 \Ricci (g)\] and the scaled Ricci flow
    \eqref{eq:RFresc} are equivalent in the following sense.\par Let
    $\left({\H^n},\ti g\left(\ti t\right)\right)_{\ti t \in [0,\ti
      T)}$ be a solution to the scaled Ricci flow.  Define
    $({\H^n},g(t))_{t \in [0,T)}$ by
    $$g(\cdot,t):= (1 + 2(n-1)t) \ti g\left( \cdot, \ti t(t)\right)  ,$$
    where $$\ti t(t):= \frac{\log( 1 + 2(n-1)t )}{2(n-1)}$$ and $T:=
    \frac{e^{2(n-1)\ti T} - 1} {2(n-1)}.$

    Then $({\H^n},g(t))_{t \in [0,T)}$ solves the Ricci flow.
    Alternatively, let $({\H^n},g(t))_{t \in [0,T)}$ be a solution to
    the Ricci flow.  Define $\left({\H^n},\ti g\left(\ti
        t\right)\right)_{\ti t \in \left[0,\ti T\right)}$ by
    $${\ti g}(\cdot, \ti t):= e^{-2(n-1)\ti t}g\left(\cdot, t\left(\ti
        t\right)\right),$$ 
    where $t(\ti t):= \frac{e^{2(n-1)\ti t} - 1} {2(n-1)}$ and $\ti
    T:= \frac{\log(1 + 2(n-1)T)}{2(n-1)}.$ Then $\tilde{g} $ solves
    the scaled Ricci flow.
  \end{lemma}
  \begin{proof}
    We prove the first claim by calculating. The second claim is shown
    in a similar way.  We calculate at $t_0$, and let $\ti t_0 :=
    \frac{\log( 1 + 2(n-1)t_0 )}{2(n-1)}$.
    \begin{align*}
      \left(\partt g\right)(\cdot,t_0) = &\, 2(n-1) \ti g\left(
        \cdot,\ti t_0 \right) +
      \left(\parttt \ti g\right)\left(\cdot, \ti t_0\right) \\
      = &\, 2(n-1) \ti g\left( \cdot,\ti t_0 \right) -2
      \Ricci\left(\ti g, \ti t_0\right)
      -2(n-1) \ti g\left( \cdot,\ti t_0 \right) \\
      = &\,-2 \Ricci\left(\ti g\left(\ti t_0\right)\right) \\
      = &\, -2 \Ricci\,(g(t_0))
    \end{align*}
    where the last line follows from the fact, that the Ricci tensor
    is invariant under scaling of the metric.
  \end{proof}

  \begin{lemma}\label{lem:real interp}
    Let $u\in C^2$ on $[0,\infty)$ or $\R$. Then \[\Vert
    Du\Vert^2_{L^\infty} \le 32\cdot \Vert u\Vert_{L^\infty} \cdot
    \left\Vert D^2u\right\Vert_{L^\infty}.\]
  \end{lemma}
  \begin{proof}
    Assume without loss of generality that $Du(0)\ge\frac12\Vert
    Du\Vert_{L^\infty} =:\frac12M$. Then $Du(x)\ge\frac14M$ for all
    $0\le x\le\frac M{4\cdot\left\Vert
        D^2u\right\Vert_{L^\infty}}$. Hence \[2\Vert u\Vert_{L^\infty}
    \ge \left|u\left(\frac M{4\cdot\left\Vert
            D^2u\right\Vert_{L^\infty}}\right) - u(0)\right| \ge\frac
    M4 \cdot \frac M{4\cdot\left\Vert D^2u\right\Vert_{L^\infty}}.\]
    The claim follows.
  \end{proof}

  \begin{lemma}\label{lemm:interp}
    Let $B$ be a compact subset of a Riemannian manifold $M$. Assume
    that $B$ has $C^2$-boundary. Let $u\in C^2(M)$. Then \[\Vert\nabla
    u\Vert^2_{L^\infty(B)} \le c(B)\cdot \Vert u\Vert_{L^\infty(B)}
    \cdot \left(\left\Vert\nabla^2u\right\Vert_{L^\infty(B)}
      +\Vert\nabla u\Vert_{L^\infty(B)}\right).\]
  \end{lemma}
  \begin{proof}
    For every point $p\in B$ and every unit vector $\xi\in T_pM$ there
    exists a curve $\gamma:[0,\infty)\to B$, parametrised by
    arc-length, such that $\gamma(0)=p$, $\gamma'(0)\in\pm\{\xi\}$ and
    \begin{align*}
      |(u\circ\gamma)'(t)| \le&\,|\nabla
      u(\gamma(t))|,\\
      |(u\circ\gamma)''(t)|
      \le&\,c(B)\cdot\left(\left|\nabla^2u(\gamma(t))\right|
        +\left|\nabla u(\gamma(t))\right|\right).
    \end{align*}
    Note that $c(B)$ depends on the curvature of $\gamma$ but can be
    chosen uniformly for all $p\in B$. Choosing $p$ and $\xi$ such
    that $\nabla u$ attains its maximum at $p$ in direction $\xi$,
    Lemma \ref{lem:real interp} yields the statement.
  \end{proof}

  \section{Euclidean Space}
  \label{eucl space app}
  We consider the situation of the main theorem, Theorem 1.3, in
  \cite{Riccistab}. Instead of a Lyapunov function involving
  $\phi_m+\psi_m-2n=\sum\limits_{i=1}^n\frac1{\lambda_i^m}
  \left(\lambda_i^m-1\right)^2$, however, we study a Lyapunov function
  involving $|g-h|^p$, $p \geq 2$. This simplifies the
  proof. \par Recall that the Ricci harmonic map heat flow with
  Euclidean background metric is
  \begin{equation*}
    \begin{split}
      \delt g_{ij}=&\ g^{ab}\nabla_a\nabla_bg_{ij}
      +\tfrac{1}{2}g^{ab}g^{pq}( \nabla_ig_{pa}\nabla_{j}g_{qb}
      +2\nabla_ag_{jp}\nabla_{q}g_{ib}\\
      &\ -2\nabla_ag_{jp}\nabla_{b}g_{iq}
      -2\nabla_jg_{pa}\nabla_{b}g_{iq}
      -2\nabla_ig_{pa}\nabla_{b}g_{jq})\ ,
    \end{split}
  \end{equation*}
  where $\nabla$ denotes covariant differentiation w.\,r.\,t.{} the
  Euclidean metric $h$. Calculating as in Lemma \ref{evolZ lem}, we
  see that
  \begin{equation}
    \label{gmh eq}
    \delt|g-h|^2-g^{ij}\nabla_i\nabla_j|g-h|^2\le
    -\left(\frac2{1+\epsilon} -9\epsilon(1+\epsilon)^2\right)|\nabla
    g|^2\le0
  \end{equation}
  if $0<\epsilon\le\frac17$. Note that there is no zeroth order term
  in the evolution equation on Euclidean space. Hence
  $|g(t)-h|\le\epsilon$ is preserved during the flow and we obtain
  long time existence, see \cite{MilesC0}. Define
  \[I_\delta^p(t):=\int\limits_{\R^n}\left(|g-h|^p-\delta\right)_+.\]
  Using \eqref{gmh eq}, and calculating as in \cite{Riccistab}, we
  get \[\frac
  d{dt}I^p_\delta(t)\le-\int\limits_{\{|g-h|^p>\delta\}}\frac
  p2\cdot\frac{2-(11+9\epsilon)(1+\epsilon)^2\epsilon}
  {1+\epsilon}\cdot|g-h|^{p-2}\cdot|\nabla g|^2\le0\] for
  $0<\epsilon\le\frac18$. The rest of the proof is similar to
  \cite{Riccistab}.  If we further restrict $p$ to $ 2 \leq p <n$ then
  we can argue as in the paper \cite{Riccistab} to prove Theorem 1.4
  of that paper.

  \section{Conformal Ricci Flow in Two Dimensions}
  \label{conf sec}
  Let us consider the Euclidean ball $B:=B_1(0)\subset\R^2$ equipped
  with the metric $(g_{ij})=\left(e^{f+u}\delta_{ij}\right)$, where
  $f=\log 4-2\log\left(1-|x|^2\right)$ and $u=u(x,t)$. For $u\equiv
  0$, we get hyperbolic space of sectional curvature equal to
  $-1$. Consider rescaled Ricci flow \[\delt
  g_{ij}=-2R_{ij}-2g_{ij}.\] As
  $R_{ij}=-\frac12\delta_{ij}\Delta_\delta(u+f)$, this is equivalent
  to
  \begin{equation}
    \label{eq:confRic}
    \dot u=e^{-u-f}\Delta_\delta u+2\left(e^{-u}-1\right)
    =e^{-u}\Delta_hu +2\left(e^{-u}-1\right) =\Delta_gu
    +2\left(e^{-u}-1\right). 
  \end{equation}
  In contrast to Theorem \ref{thm:Ricci}, we do not have to assume
  that the eigenvalues $(\lambda_i)$ of $g(0)$ with respect to $h$ are
  close to one. This is similar to \cite[Theorem
  A.1]{Riccistab}. There, however, we had to assume that
  $\lambda_i(x,0)\to1$ for $|x|\to\infty$ in order to obtain
  convergence to $\R^2$, see \cite[Theorem A.2]{Riccistab}.
  \begin{theorem}
    Let $u_0\in C^0(B)$ satisfy $\Vert
    u_0\Vert_{L^\infty}<\infty$. Then there exists $u\in
    C^\infty(B\times(0,\infty))$ solving \eqref{eq:confRic} such that
    $u(\cdot,t)\to u_0$ in $C^0_{\loc}(B)$ as $t\searrow0$. As
    $t\to\infty$, $u(\cdot,t)\to 0$ exponentially in $C^\infty$
    w.\,r.\,t.{} the hyperbolic metric. 
  \end{theorem}
  If an arbitrary solution $u$ is uniformly bounded for
  small times, we also get exponential convergence. 
  \begin{proof} 
    Assume $|u_0| \leq c_0$.  Mollify and modify $u_0$ to $u^i_0$ with
    $|u^i_0| \leq 2 c_0$ and $u^i_0 = 0$ near $\boundary
    B_{1-\frac1i}(0)$ and $u^i_0 = u_0 $ on $B_{1-\frac2i}(0)$.  We
    can construct solutions $u^i: B_{1-\frac1i}(0) \times [0, T_i) $
    to \ref{eq:confRic} with $u^i(\cdot,0) = u_0 $ on $B_{1-\frac 1
      i}(0)$ and $u_i(\cdot,t) = 0$ on $\partial B_{1-\frac 1 i}(0)$
    using the arguments presented in Chapter VI of
    \cite{LadySoloUral}.  These solutions remain bounded by $2c_0$
    from the maximum principle.  Hence, the arguments of Chapter VI of
    \cite{LadySoloUral} imply that $T_i = \infty$.
    
    Spatially constant barriers $b=b(t)=\log\left(1+ae^{-2t}\right)$,
    $a>-1$, solving \eqref{eq:confRic} converge exponentially to $0$
    as $t\to\infty$. Hence the maximum principle applied to each $u^i$
    on $B_{1-\frac1i}(0)$ implies that the $u_i$ remain
    uniformly bounded and go exponentially to zero (uniformly in $i$)
    as $t \to \infty$.
    
    Now we address smooth convergence: Writing $l^i := e^{u^i}$ we
    obtain the evolution equation
    \[\partt l^i= (1/l^i) \Delta_h l^i 
    - (1/(l^i)^2)|\gradh l^i|^2 + 2\left( 1 - l^i \right).\] We can
    assume without loss of generality, that $|l_i-1| \leq \ep$ for
    some small $\ep$. The interior estimates of Lemma \ref{lem:intest}
    hold here as the equation for $l^i$ has the same form as the
    equation studied in Theorem 4.3 of \cite{MilesC0}.  Hence, by
    taking a diagonal subsequence, we get a solution $l = e^u$ which
    approaches $1$ exponentially.  Interpolating between the
    $C^0$-norm and $C^k$-norms and using Lemma \ref{lem:intest} again,
    we see that $l$ approaches $1$ in all $C^k$-norms exponentially .
  \end{proof}

  To treat the question of uniqueness of such solutions we work in the
  unrescaled setting. Note that by Lemma \ref{lem:equiv} this is
  equivalent to the rescaled equation.\\
  With respect to the hyperbolic metric $h$ on $\mathbb{H}^2$ as
  a background metric, a solution $e^{u(p,t)}h$ to the Ricci flow
  satisfies
  \begin{equation}\label{eq:conformunresc}
    \dot{u} = e^{-u} \Delta_h u + 2 e^{-u}\ .
  \end{equation}
  We first prove a noncompact maximum principle.
  \begin{lemma}
    Let $v\in C^\infty\left(\mathbb{H}^2, [0,T)\right)$ be a bounded
    solution to
    \begin{equation}\label{eq:equinterpolated}
      \dot{v}\leq a \Delta_h v + c\, v
    \end{equation}
    with $a,c \in L^\infty\left(\mathbb{H}^2,[0,T)\right)$, $a>0$. If
    $v(\cdot,0)\leq 0$ then $v(\cdot, t)\leq 0$ for all $t\in [0,T)$.
  \end{lemma}
  \begin{proof}
    Pick a fixed point $p_0\in \mathbb{H}^2$ and let $r(\cdot):=
    \text{dist}_h(\cdot, p_0)$. Then the function $\rho:=
    \sqrt{r^2+1}$ is a smooth function on $\mathbb{H}^2$ with $\rho(p)
    \rightarrow \infty$ as $p\rightarrow \infty$ and
    \begin{equation*}
      a\Delta_h\rho \leq C
    \end{equation*}
    for a constant $C>0$. Let us first assume that $v$ satisfies
    $$\dot{v}\leq a \Delta_h v - c' v$$
    with a function $c'\geq 0$. Then for any $\delta>0$ the function
    $$w:= v-\delta \rho- 2\delta C t -\delta $$
    satisfies at the first non-negative interior maximum
    $$\dot{w} < a \Delta_h w .$$
    Since $w(\cdot,t) \rightarrow -\infty$ as $p\rightarrow \infty$ an
    application of the maximum principle proves the estimate as
    $\delta \rightarrow 0$. In the
    general case let $|c(p,t)|\leq K$ and $v':=e^{-Kt}v$ wich satisfies
    $$\dot{v'} \leq a \Delta_h v' -(K-c)v' \ .$$
    The previous estimate can be applied.
  \end{proof}
  This gives us a uniqueness statement.
  \begin{lemma}
    Let $u,\tilde{u}\in C^\infty(\mathbb{H}^2,(0,T))\cap
    C^0(\mathbb{H}^2,[0,T))$ be two bounded solutions of
    \eqref{eq:conformunresc} s.t.
    $$ u(\cdot ,t) \rightarrow  u_0 \ \text{and} \ \tilde{u}(\cdot ,t) 
    \rightarrow  u_0$$
    uniformly as $t\rightarrow 0$ for some continuous function $u_0$ on
    $\mathbb{H}^2$. Then $u\equiv \tilde{u}$.
  \end{lemma}
  \begin{proof}
    Define for $\gamma>0$
    $$ u_\gamma(p,t):= u(p,e^{-\gamma}t)+\gamma\ . $$
    Then $u_\gamma$ again solves \eqref{eq:conformunresc} with inital
    values $u_0+\gamma$. Since the initial values are attained
    uniformly we have $u_\gamma > \tilde{u}$ for a short time interval
    $[0,2\delta], \delta >0$. By interior estimates as in Lemma
    \ref{lem:intest} the functions $u_\gamma,\tilde{u}$ are bounded
    uniformly in $C^\infty$ on time intervals
    $[\delta,T)$. Interpolating between the two solutions, we see that
    the difference satisfies an equation of the form
    \eqref{eq:equinterpolated}, to which the noncompact maximum
    principle applies. Thus $u_\gamma > \tilde{u}$ for all $\gamma >0$
    and $\gamma \rightarrow 0$ gives the desired estimate.
  \end{proof}

  If $u_0$ is unifomly continous then also $u(\cdot,t)$ converges
  uniformly as $t\rightarrow 0$.
  \begin{lemma}
    Let $u \in C^\infty\left(\mathbb{H}^2,(0,T)\right)\cap
    C^0\left(\mathbb{H}^2,[0,T)\right)$ be a bounded solution of
    \eqref{eq:conformunresc} s.t. $u_0:=u(\cdot,0)$ is uniformly
    continuous. Then $u(\cdot, t)\rightarrow u_0$ uniformly.
  \end{lemma}
  \begin{proof}
    Pick a point $p_0 \in \mathbb{H}^2$ and let $u_{p_0}(t)$ be the solution
    to \eqref{eq:conformunresc}, which is constant in space and has
    initial value $u(p_0)$. It is a direct computation to check that
    \begin{equation*}
      \begin{split}
        \delt (u-u_{p_0})^2 =&\ e^{-u}\Delta_h(u-u_{p_0})^2 -
        2e^{-u}|\nabla^hu|^2 + 2(u-u_{p_0})(e^{-u}-e^{-u_{p_0}})\\
        \leq &\ e^{-u}\Delta_h(u-u_{p_0})^2 -
        2e^{-u}|\nabla^hu|^2 .
      \end{split}
    \end{equation*}
    By a similar argument as in Lemma \ref{lem:intcloseness}, but now
    for small radii, one obtains the desired closeness-estimate. 
  \end{proof}
  \begin{corollary}
    Any bounded solution to \eqref{eq:confRic} in
    $C^\infty\left(\mathbb{H}^2,(0,\infty)\right)\cap
    C^0\left(\mathbb{H}^2,[0,T)\right)$ with initial data $u_0$, which
    is uniformly continuous with respect to the hyperbolic metric, is
    unique. As $t\rightarrow \infty, u(\cdot,t)\rightarrow 0$
    exponentially in $C^\infty$.
  \end{corollary}

\end{appendix}

\bibliographystyle{amsplain}

\begin{thebibliography}{10}

\bibitem{ChLuNi} Bennett Chow, Peng Lu, and Lei Ni, \emph{Hamilton's
    {R}icci flow}, Graduate Studies in Mathematics, vol.~77, American
  Mathematical Society, Providence, RI, 2006.

\bibitem{DeTurck} Dennis~M. DeTurck, \emph{Deforming metrics in the
    direction of their {R}icci tensors}, J. Differential
  Geom. \textbf{18} (1983), no.~1, 157--162.

\bibitem{HamiltonThree} Richard~S. Hamilton, \emph{Three-manifolds
    with positive {R}icci curvature}, J.  Differential
  Geom. \textbf{17} (1982), no.~2, 255--306.

\bibitem{LadySoloUral} Olga~A. Lady{\v{z}}enskaja,
  Vsevolod~A. Solonnikov, and Nina~N.  Ural{\cprime}ceva, \emph{Linear
    and quasilinear equations of parabolic type}, Translated from the
  Russian by S. Smith. Translations of Mathematical Monographs,
  Vol. 23, American Mathematical Society, Providence, R.I., 1967.

\bibitem{LiYinHypRic} Haozhao Li and Hao Yin, \emph{On stability of
    the hyperbolic space form under the normalized {R}icci flow},
  Int. Math. Res. Not. IMRN (2010), doi:10.1093/imrn/rnp232.

\bibitem{McKean} H.~P. McKean, \emph{An upper bound to the spectrum of
    {$\Delta $} on a manifold of negative curvature}, J. Differential
  Geometry \textbf{4} (1970), 359--366.

\bibitem{Riccistab} Oliver~C. Schn{\"u}rer, Felix Schulze, and Miles
  Simon, \emph{Stability of {E}uclidean space under {R}icci flow},
  Comm. Anal. Geom. \textbf{16} (2008), no.~1, 127--158.

\bibitem{ShiJDG1989} Wan-Xiong Shi, \emph{Deforming the metric on
    complete {R}iemannian manifolds}, J. Differential
  Geom. \textbf{30} (1989), no.~1, 223--301.

\bibitem{MilesC0} Miles Simon, \emph{Deformation of {$C\sp 0$}
    {R}iemannian metrics in the direction of their {R}icci curvature},
  Comm. Anal. Geom. \textbf{10} (2002), no.~5, 1033--1074.

\bibitem{Suneeta} V.~Suneeta, \emph{Investigating the off-shell
    stability of anti-de {S}itter space in string theory}, Classical
  Quantum Gravity \textbf{26} (2009), no.~3, 035023, 18.

\bibitem{YeEinstein} Rugang Ye, \emph{Ricci flow, {E}instein metrics
    and space forms}, Trans. Amer.  Math. Soc. \textbf{338} (1993),
  no.~2, 871--896.

\end{thebibliography}
\def\cprime{$'$}
\providecommand{\bysame}{\leavevmode\hbox to3em{\hrulefill}\thinspace}
\providecommand{\MR}{\relax\ifhmode\unskip\space\fi MR }
\providecommand{\MRhref}[2]{%
  \href{http://www.ams.org/mathscinet-getitem?mr=#1}{#2}
}
\providecommand{\href}[2]{#2}

\end{document}